% This is a list of corrections for the book: J. Noguchi and T. Ochiai,
% Geometric Function Theory in Several Complex Variables,
% xi + 282 pp., Math.\ Monographs Vol.\ {\bf 80},
% Amer.\ Math.\ Soc., Providence, 1990.
% 
% The authors hope that this distribution will be helpful for readers
% to avoid unnecessary confusions.
% 
%%%%%%%%%%%%%%%%%%%% This is a Plain-Tex file.  %%%%%%%%%%%%%%%%%%%
%**start of header
% CALL 12pt FONTS
\font\twelverm=cmr12
\font\twelvebf=cmbx12
\font\twelveit=cmti12
\def\bigfont{		% DEFINE 12pt FONTS
	\let\rm=\twelverm
	\let\bf=\twelvebf
	\let\it=\twelveit
	\rm}
% for Letter size
\vsize=22true cm
\hsize=17true cm
% for A4 size
%\vsize=24true cm
%\hsize=16true cm
\headline={\ifnum\pageno=1 \hfil \else\hfil {$\folio$} \hfil\fi}
\baselineskip=10pt
\parskip=0pt
%**end of header

\centerline{\bigfont\bf Corrections}
\centerline{\bf (Version, 31st March 1996)}
\centerline{\bigfont\bf Geometric Function Theory in Several Complex Variables}
\medskip
\centerline{\bf J. Noguchi and T. Ochiai}
%\medskip
%\centerline{\bf Version, 31st March 1996}
\medskip
\noindent p. xi, $\downarrow 3$: \S1 $\Longrightarrow$ \S2\hfil\break
p. 14, $\uparrow 3$: $\tilde \gamma_i \Longrightarrow {\dot {\tilde \gamma}}_i$\hfil\break
p. 15, $\downarrow 2,5,6$: $\tilde \gamma_i \Longrightarrow {\dot {\tilde \gamma}}_i$\hfil\break
p. 17, $\downarrow 16$: $F \Longrightarrow G$\hfil\break
p. 19, $\uparrow 6$: 0 Then $\Longrightarrow$ 0. Then\hfil\break
p. 30, $\uparrow 1$: $F \Longrightarrow L$\hfil\break
p. 31, $\downarrow 15$: $B(3/4) \Longrightarrow \overline{B(3/4)}$\hfil\break
p. 32, $\downarrow 3$: $\Delta^* \Longrightarrow \Delta^* (1)$\hfil\break
p. 32, $\downarrow 6$: $\empty^\prime \quad(\hbox{2 places}) \Longrightarrow$ delete\hfil\break
p. 32, $\downarrow 13$: $\min \Longrightarrow$ delete\hfil\break
p. 34, $\downarrow 5$: $h(z) \Longrightarrow |h(z)|$ (2 places)\hfil\break
p. 39, $\uparrow 11$: (1.6.6) $\Longrightarrow$ (1.6.5)\hfil\break
p. 41, $\downarrow 11$: $U_\nu \Longrightarrow V_\nu$\hfil\break
p. 44, $\downarrow 11$: $\psi_{i} \colon (z_{i}^{0},...,
\buildrel {\wedge} \over {z_i^i} ,..., z_{i}^{m}) \in {\bf C}^{m}
\rightarrow \rho (z_{i}^{0},...,
\buildrel {i-{\rm th}} \over 1 ,..., z^{m}) \in U_{i}$
$\Longrightarrow$
\hfil\break
\hglue60pt
$\psi_{i} \colon \rho (z_{i}^{1},...,
\buildrel {i-{\rm th}} \over 1 ,..., z_{i}^{m}) \in U_{i} \rightarrow
(z_{i}^{1},...,
\buildrel \wedge \over {z_{i}^{i}} ,..., z_{i}^{m}) \in {\bf C}^{m-1}$\hfil\break
p. 44, $\downarrow 12$: ${\bf C}^{m}$ $\Longrightarrow$ $E$\hfil\break
p. 44, $\downarrow 14$: ${\bf C}^m$  $\Longrightarrow$  ${\bf C}^{m-1}$\hfil\break
p. 44, $\downarrow 15$: ${\bf C}^{m}$ $\Longrightarrow$ $E$\hfil\break
p. 44, $\downarrow 18$: $z^0_i$  $\Longrightarrow$  $z^1_i$\hfil\break
p. 60, $\downarrow 13$: ${\bf 1}_{M}$ $\Longrightarrow$ ${\bf 1}_{U}$\hfil\break
p. 64, $\downarrow 20$: $H$ $\Longrightarrow$ $H_{i}$\hfil\break
p. 65, $\uparrow 10$: $->$ $\Longrightarrow$ $\rightarrow$\hfil\break
p. 68, $\downarrow 5$: $X$ $\Longrightarrow$ $x$\hfil\break
p. 70, $\downarrow 11$: $= m$ $\Longrightarrow$ $= 2m$\hfil\break
p. 71, $\uparrow 8$: ${\bf K}_{{\bf C}^{m}}$ $\Longrightarrow$
${\bf K}({\bf C}^{m})$\hfil\break
p. 72, $\uparrow 11$: $)^{k-m-1}$ $\Longrightarrow$ $)^{d-m-1}$\hfil\break
p. 72, $\uparrow 11$: $k$ $\Longrightarrow$ $d$\hfil\break
p. 72, $\uparrow 11$: $\{$ $\Longrightarrow$ $[$\hfil\break
p. 72, $\uparrow 11$: $\}$ $\Longrightarrow$ $]$\hfil\break
p. 73, $\uparrow 3$: $\Omega \geq 0$) for $\Longrightarrow$
$\Omega (x) \geq 0$) for\hfil\break
p. 76, $\downarrow 4$: $tr_{m})$ $\Longrightarrow$ $tr_{m})(z)$\hfil\break
p. 76, $\downarrow 15$: $- \partial \bar\partial$ $\Longrightarrow$ $-i \partial \bar\partial$\hfil\break
p. 76, $\uparrow 9$: $\Omega^{m}$ $\Longrightarrow$ $\Omega$\hfil\break
p. 78, $\downarrow 13$: $a(y) < b(y)$ $\Longrightarrow$ $a(y) \leq b(y)$\hfil\break
p. 79, $\downarrow 8$: $\Lambda^m$ $\Longrightarrow$ $\Lambda$\hfil\break
p. 79, $\downarrow 12$: $j=1$ $\Longrightarrow$ $j=2$\hfil\break
p. 81, $\downarrow 2$: $\int_M$ $\Longrightarrow$ $\int_{B'}$\hfil\break
p. 81, $\uparrow 9$: $\displaystyle \sum_{j=1}^\infty \Psi_M$ $\Longrightarrow$ $\displaystyle \sum_{j=1}^\infty \int_{f_j(E_j)} \Psi_M$\hfil\break
p. 83, $\downarrow 16$: (5.4.1) $\Longrightarrow$ (5.5.1)\hfil\break
p. 83, $\downarrow 21$: (2.1.8) $\Longrightarrow$ (2.1.22)\hfil\break
p. 83, $\downarrow 22$: (2.1.9) $\Longrightarrow$ (2.1.23)\hfil\break
p. 84, $\downarrow 10$: Put equation number (2.6.3).\hfil\break
p. 86, $\downarrow 15$: (iv) $\Longrightarrow$ (d)\hfil\break
p. 86, $\downarrow 16$: (3.3.43) $\Longrightarrow$ (3.3.44)\hfil\break
p. 87, $\downarrow 14$: $\Phi$ $\Longrightarrow$ $\Psi$\hfil\break
p. 88, $\uparrow 9$: $\| f^{\ast} \xi \|$ $\Longrightarrow$
$\| f^{\ast} \xi \|_{N}$\hfil\break
p. 88, $\uparrow 6$: $\| \xi \|$ $\Longrightarrow$
$\| \xi \|_{M}$\hfil\break
p. 91, $\downarrow 7$: {\it manifolds} $\Longrightarrow$ {\it spaces}\hfil\break
p. 91, $\downarrow 12$: of Lang $\Longrightarrow$ (delete)\hfil\break
p. 91, $\downarrow 12$: K\"ahler $\Longrightarrow$ a K\"ahler manifold\hfil\break
p. 91, $\downarrow 13$: compact $\Longrightarrow$ (delete)\hfil\break
p. 91, $\downarrow 13$: Moisezon $\Longrightarrow$ K\"ahler manifold\hfil\break
p. 99, $\uparrow 1$: $\| \phi \|_0 T(\phi).$ $\Longrightarrow$ $\| \phi \|_0 T(\phi_A) \pm T(\phi)$.\hfil\break
p. 111, $\downarrow 4$: {\bf real current.} $\Longrightarrow$ {\bf real current}, and $p=q$.\hfil\break
p. 113, $\downarrow 7$: (3.2.14) $\Longrightarrow$ (3.1.14)\hfil\break
p. 114, $\downarrow 4$: positive distributions $\Longrightarrow$ distributions
of order 0
p. 114, $\downarrow 13$: $\sigma$ $\Longrightarrow$ $\sigma_k$\hfil\break
p. 117, $\uparrow 1$: $\displaystyle\int_{\overline {B(r_{2}) - B(r_{1})}}$
$\Longrightarrow$ $\displaystyle\int_{B(r_{2}) - \overline {B(r_{1})}}$
\hfil\break
p. 118, $\uparrow 5$: $< T$ $\Longrightarrow$ $\leq T$\hfil\break
p. 120, $\downarrow 11$: (two places) $1 \over {r^{2k}}$ $\Longrightarrow$
(1'st) $1 \over {r_2^{2k}}$ ; (2'nd) $1 \over {r_1^{2k}}$\hfil\break
p. 120, $\uparrow 8$: $<$ $\Longrightarrow$ $\leq$\hfil\break
p. 121, $\uparrow 10$: $\epsilon^{i\theta}$ $\Longrightarrow$ $\epsilon e^{i\theta}$\hfil\break
p. 122, $\downarrow^{} 12$: {\it increasing} $\Longrightarrow$ {\it decreasing}
\hfil\break
p. 122, $\downarrow 13$: $u_{j}(z) \leq u_{j+1}$ $\Longrightarrow$
$u_{j}(z) \geq u_{j+1}(z)$\hfil\break
p. 122, $\downarrow 13$: $u_{j}$ $\Longrightarrow$ $u_{j}(z)$\hfil\break
p. 122, $\uparrow 8$: above $\Longrightarrow$ above on $K$\hfil\break
p. 128, $\uparrow 4$: $[u]$ $\Longrightarrow$ $[u_\epsilon ]$\hfil\break
p. 129, $\downarrow 10, 12$: (two places) $\leq$ $\Longrightarrow$ $=$\hfil\break
p. 130, $\uparrow 8$: $\Delta(1)$ $\Longrightarrow$ $\Delta (R)$\hfil\break
p. 140, $\downarrow 11$: $b_{U}$ $\Longrightarrow$ $b_{V}$\hfil\break
p. 144, $\downarrow 15$: $\hbox{\bigfont \it I}_{M,x}$ $\Longrightarrow$
$\hbox{\bigfont\it I}_{X,x}$\hfil\break
p. 145, $\downarrow 13$: (ii) $\Longrightarrow$ (iii)\hfil\break
p. 148, $\downarrow 12$: $=m($ $\Longrightarrow$ $=-m($\hfil\break
p. 149, $\uparrow 13$: $f$ $\Longrightarrow$ $F$\hfil\break
p. 150, $\downarrow 3$: $a$ $\Longrightarrow$ $b$\hfil\break
p. 150, $\uparrow 8$: $e^{-u_{i}}$ $\Longrightarrow$ $e^{-u_{(n)i}}$\hfil\break
p. 151, $\downarrow 7$: $f^{-1}(x)$ $\Longrightarrow$ $f^{-1}(f(x))$\hfil\break
p. 152, $\uparrow 13$: $M$ $\Longrightarrow$ $X$\hfil\break
p. 154, $\uparrow 3$: $U \cap U_2 = \emptyset$ $\Longrightarrow$
$U_1 \cap U_2 = \emptyset$\hfil\break
p. 160, $\downarrow 4$: $\tilde {\bf C}^{{n+1}^{n+1}}$ $\Longrightarrow$ 
$\tilde {\bf C}^{n+1}$\hfil\break
p. 166, $\downarrow 6$: $V$ $\Longrightarrow$ $A$\hfil\break
p. 188, $\downarrow 7$: $B(r)$ $\Longrightarrow$ $\overline {B(r)}$\hfil\break
p. 206, $\downarrow 2$: (5.4.17) $\Longrightarrow$ (5.5.17)\hfil\break
p. 207, $\uparrow 1$: Put equation number (5.5.22).\hfil\break
p. 214, $\uparrow 6$: (5.4.6) $\Longrightarrow$ (5.5.46)\hfil\break
p. 242, $\downarrow 6$: $X$ $\Longrightarrow$ $Y$\hfil\break
p. 244, $\uparrow 7$: (6.2.5) $\Longrightarrow$ (6.1.5)\hfil\break
p. 247, $\downarrow 6$: states $\Longrightarrow$ stated\hfil\break
p. 247, $\uparrow 3$: variety. $\Longrightarrow$ variety $A.$\hfil\break
p. 247, $\uparrow 2$: {\it any} $\Longrightarrow$ {\it Any}\hfil\break
p. 254, $\uparrow 11$: $u(x)$ $\Longrightarrow$ $|u(x)|$\hfil\break
p. 254, $\uparrow 7$: $u_{r}$ $\Longrightarrow$ ${\buildrel \wedge \over u}_{r}$\hfil\break
p. 254, $\uparrow 5$: $u^{m}$ $\Longrightarrow$ $y^{m}$\hfil\break
p. 255, $\downarrow 2 \sim 4$: $n$ $\Longrightarrow$ $m$\hfil\break
p. 256, $\downarrow 2$: $x$ $\Longrightarrow$ $z$\hfil\break
p. 258, $\uparrow 12$: $dz^{j} \wedge d \bar z ^{k}$ $\Longrightarrow$
$T_{j \bar k} dz^{j} \wedge d \bar z ^{k}$\hfil\break
p. 260, $\uparrow 6$: $q(q; a, z)$ $\Longrightarrow$ $e(q; a, z)$
\hfil\break
p. 261, $\downarrow 9 \sim 10$: $q,$ $\Longrightarrow$ $q;$\hfil\break
p. 261, $\downarrow 11$: $(($ $\Longrightarrow$ $($\hfil\break
p. 263, $\uparrow 3$: $r \rightarrow \infty$ $\Longrightarrow$
$\nu \rightarrow \infty$\hfil\break
p. 264, $\uparrow 2 \sim 1$: By $\sim$ (iii), $\Longrightarrow$
By (5.2.25), Corollary (5.2.30), condition (iii) and Lemma (3.2),\hfil\break
p. 265, $\downarrow 1 \sim 2$: Change these 2 lines by the following 2 lines:
$$T(r_{\nu}, A^{\pm 1} ) \leq
T(r_{\nu}, F) + T(r_{\nu} , G) + O(1) = o(r_{\nu}^{q+1}).$$
Thus Theorem (5.3.13) implies that
$\log M(r_{\nu}, A^{\pm 1}) = o(r_{\nu}^{q+1}).$\hfil\break
p. 265, $\downarrow 4$: $|{\rm Re} \, \log A(z) |$ $\Longrightarrow$
${\rm sup} \{|{\rm Re} \, \log A(z)|; \; z \in B(r_{\nu})\}$\hfil\break
p. 269, $\uparrow 12 \sim 11$: ``Compact $\sim$ (1987). $\Longrightarrow$
``A finiteness criterion for compact varieties of
surjective holomorphic mappings,'' Kodai Math. J. {\bf 13},
pp. 373-376 (1990).\hfil\break
p. 269, $\uparrow 8$: K. $\Longrightarrow$ Y.\hfil\break

\bye